\documentclass[a4paper]{article}      
\usepackage{amsmath,amssymb,amsfonts} 
\usepackage{graphics}                 
\usepackage{color}                    
\usepackage{hyperref}                 
\newtheorem{theorem}{Theorem}[section]
\newtheorem{corollary}{Corollary}[section]

\newtheorem{example}{Example}[section]

\newtheorem{proposition}{Proposition}

\newcommand{\qed}{\hfill\rule{2mm}{2mm}}  
\newenvironment{proof}{
\begin{trivlist}
\item[\hspace{\labelsep}{\bf\noindent Proof. }]
}{\qed\end{trivlist}}
\title{\bf A symbolic method for $k$-statistics}
\author{E. Di Nardo,
D. Senato\footnote{Dipartimento di Matematica,
 Universit\`a degli Studi della Basilicata
 C.da Macchia Romana, I-85100 Potenza,
  E-mail: \tt \{dinardo,senato\}@unibas.it}}
\date{}
\begin{document}
\setlength{\baselineskip}{16pt}
\maketitle
\begin{abstract}
Trough the classical umbral calculus, we provide new, compact and
easy to handle expressions of $k$-statistics, and more in general
of $U$-statistics. In addition such a symbolic method can be
naturally extended to multivariate case and to generalized
$k$-statistics.
\end{abstract}
%
\section{Introduction}
In 1929, Fisher \cite{Fisher} introduced the $k$-statistics as new
symmetric functions of the random sample. The aim of Fisher was to
estimate the cumulants by free-distribution methods without using
moment estimators. He used only combinatorial techniques.
$k$-statistics are related to the power sum symmetric polynomials
whose variables are the random variables (r.v.'s) of the sample.
The Fisher point of view is described with a wealth of details by
Kendall and Stuart in \cite{Stuart-Ord}. The method is
straightforward enough,  however his execution leads to some
intricate computations and some cumbersome expressions, except in
very simple cases. This is why many authors tried  to
simplify the matter later.
\par
One of the most relevant contributions has been given by Speed
\cite{Speed2} in the Eighty. Speed resumed the Doubilet approach
to symmetric functions \cite{Doubilet}, exploiting symmetric
functions labelled by partition of a set rather than by partition
of an integer. The entries in the transition matrices are computed
via M{\oe}bius function generalizing  and simplifying the
presentation of $k$-statistics theory. Nevertheless, in order to
extend such theory to generalized $k-$ statistics,  Speed must
resort the tensor approach introduced by Kaplan in 1952.
\par
We follow a different point of view by using the high
computational potential of the classical umbral calculus.
This symbolic method was introduced by Rota and Taylor in 1994 \cite{SIAM} and further
developed in \cite{Dinardo} and \cite{Dinardo1}. From a
combinatorial perspective, we revisit the Fisher theory as exposed by
Kendall and Stuart, taking into account the Doubilet approach to
symmetric functions. The umbral calculus offers a nimble syntax
method that allows both the computation without using the
M{\oe}bius function and a natural extension to the multivariate
case without bringing the tensor device into. What is more, this
language clarifies the role played by the power sum symmetric
polynomials in the expressions of $k$-statistics.
\par
After recalling in Section 2 the strictly necessary umbral
background, in Section 3 we put in umbral setting the four
classical bases of symmetric polynomial algebra. Then we proceed
in defining a general procedure to write down $U$-statistics. Such
a procedure is given in full details in Section 5, where analogous
formulae for $k$-statistics, $h$-statistics and multivariate
$k$-statistics are given.
\section{The umbral calculus language}
We start presenting the formal setting of the umbral calculus as
introduced by Rota and Taylor in \cite{SIAM} and further developed
in \cite{Dinardo} and \cite{Dinardo1}. We shall confine our
exposition to what is necessary to the aims of this paper.
\par
The umbral calculus is a syntax consisting of the following data:
an alphabet $A=\{\alpha,\beta, \ldots \}$ whose elements are named
{\it umbrae}; a commutative integral domain $R$ whose quotient
field is of characteristic zero; a linear functional $E,$ called
{\it evaluation}, defined on the polynomial ring $R[A]$ and taking
values in $R$ and such that $E[1]=1,$ $E[\alpha^i \beta^j \cdots
\gamma^k] = E[\alpha^i]E[\beta^j] \cdots E[\gamma^k]$ for any set
of distinct umbrae in $A$ and for $i,j,\ldots,k$ non-negative
integers ({\it uncorrelation property}); an element $\epsilon \in
A,$ called {\it augmentation}, such that $E[\epsilon^n] = 0$ for
$n \geq 1;$ an element $u \in A,$ called {\it unity} umbra, such
that $E[u^n]=1,$ for $n \geq 1.$
\par
The {\it support} of an {\it umbral polynomial} $p \in R[A]$
is the set of all umbrae occurring in $p.$ Two umbral polynomials
are said to be {\it uncorrelated} when their supports are disjoint.
\par
An umbra carries the structure of a r.v., while making no
reference to a probability space if the evaluation $E$ is
considered as the expectation operator. The {\it moments} of an
umbra $\alpha$ are the elements $a_n \in R$ such that
$E[\alpha^n]=a_n$ for $n \geq 1.$ We said that the umbra $\alpha$
{\it represents} the sequence $1,a_1,a_2,\ldots.$ The {\it
singleton umbra} $\chi$ is the umbra whose moments are all zero,
but the first equal to $1.$ The {\it factorial moments} of an
umbra $\alpha$ are the elements $a_{(n)} \in R$ corresponding to
the umbral polynomials $(\alpha)_n = \alpha (\alpha -1) \cdots
(\alpha-n+1), n \geq 1$ via the evaluation $E,$ i.e.
$E[(\alpha)_n]=a_{(n)}.$ The {\it Bell umbra} $\beta$ is the umbra
such that $E[(\beta)_n]=1, n \geq 1.$ The umbra $\beta$ and the
umbra $\chi$ allow respectively a symbolic tool in order to
handling composition and inversion of formal power series (for a
detailed exposition of such two umbrae with their properties see
\cite{Dinardo} and \cite{Dinardo1}).
\par
Two umbrae  $\alpha$ and $\gamma$ are said to be {\it similar}
when $E[\alpha^n]=E[\gamma^n]$ for all $n \geq 1,$ and we will
write $\alpha \equiv \gamma.$ Then two similar umbrae represent
the same moment sequence. Given a sequence $1,a_1,a_2,\ldots$ in
$R$ this is represented by infinitely many distinct and thus
similar umbrae. Such a circumstance permits to deal with moment
products by means of the uncorrelation property. Indeed, it is
$a_i a_j \ne E[\alpha^i \alpha^j]$ with $a_i = E[\alpha^i]$ and
$a_j=E[\alpha^j],$ as well as $a_i a_j = E[\alpha^i \alpha^{'j}]$
with $\alpha \equiv \alpha^{'}$ and $\alpha^{'}$ uncorrelated with
$\alpha.$ So, given $n \in N,$ the sequence
$$\sum_{k=0}^n \left( \begin{array}{c}
n \\
k \end{array} \right) a_{n-k} a_k$$
gives the moments of $\alpha+\alpha^{'}.$
\par
Let $p$ and $q$ be two umbral polynomials. We said that $p$ is
{\it umbrally equivalent} to $q$ iff $E[p]=E[q],$ in symbols $p
\simeq q.$ This last equivalence relation turns out to be very
useful in defining and handling umbra generating functions. We
point out that all operations among umbrae correspond to analogous
operation in the algebra of generating functions. Nevertheless, in
the following we make no mention of generating functions
(referring \cite{Taylor1} to a formal exposition).
\par
Thanks to the introduction of similarity notion, it is possible
to extend the alphabet $A$ with the so-called {\it auxiliary umbrae}
derived from operations among similar umbrae. This leads to the
construction of a {\it saturated umbral calculus} in which the auxiliary
umbrae are treated as elements of the alphabet (cf. \cite{SIAM}).
Let $\alpha^{'},\alpha^{''},\ldots, \alpha^{'''}$ be $n$
uncorrelated umbrae similar to an umbra $\alpha.$ The symbol $n.\alpha$ denotes
an auxiliary umbra similar to the sum $\alpha^{'}+\alpha^{''}+\cdots+\alpha^{'''}.$
The symbol $\alpha^{.n}$ denotes an auxiliary umbra similar to
the product $\alpha^{'}\alpha^{''} \cdots \alpha^{'''}.$
Properties of such auxiliary umbrae are extensively described
in \cite{Dinardo} and they will be recalled whenever it is
necessary. We will assume disjoint both the support of
$n.\alpha, m.\alpha$ and $\alpha^{.n},\alpha^{.m}$
whenever $n \ne m.$ If $p$ and $q$ are correlated umbral
polynomials, then $n.p \simeq p_1+\cdots+p_n$ is
correlated to $n.q \simeq q_1+\cdots+q_n,$ and
$p_i$ is correlated to $q_i$ but uncorrelated to $q_j$
with $i \ne j.$ In \cite{Dinardo}, the following identity
is stated:
\begin{equation}
E[(n.\alpha)^i]=\sum_{k=1}^i (n)_k B_{i,k}(a_1,a_2,\ldots,a_{i-k+1})
\quad i \geq 1,
\label{(momdot)}
\end{equation}
where $B_{i,k}$ are the (incomplete) exponential Bell polynomials
and $a_i$ is the $i$-th moment of $\alpha.$ Moreover , it is easy
to verify that $E[(\alpha^{.n})^i]=a_i^n$ for $i \geq 0.$
\par
Two umbrae $\alpha$ and $\gamma$ are said to be {\it inverse} to
each other when $\alpha+\gamma \equiv \varepsilon.$ The
inverse of the umbra $\alpha$ is denoted by $-1.\alpha.$ Note that, in dealing with a
saturated umbral calculus, the inverse of an umbra is not unique,
but any two inverse umbrae of the same umbra are similar.
\par
Replacing the integer $n$ in $n.\alpha$ with an umbra $\gamma,$
we obtain the auxiliary umbra $\gamma.\alpha$ whose moments are
\begin{equation}
E[(\gamma.\alpha)^i]=\sum_{k=1}^i g_{(k)} B_{i,k}(a_1,a_2,\ldots,a_{i-k+1})
\quad i \geq 1,
\label{(ombdot)}
\end{equation}
where $g_{(k)}$ are the factorial moments of $\gamma.$ In
particular $\beta.\alpha$ is called {\it $\alpha-$partition
umbra} and its moments are the (complete) exponential Bell
polynomials (cf. \cite{Dinardo}). Moreover $\chi.\alpha$ is
called {\it $\alpha-$cumulant umbra} and $\alpha.\chi$ is called
{\it $\alpha-$factorial umbra}, with moments equal to the
factorial moments of $\alpha$ (cf. \cite{Dinardo1}). In particular
it is
\begin{equation}
\beta.\chi \equiv u \equiv \chi.\beta.
\label{(chibeta)}
\end{equation}
Again, replacing the umbra $\gamma$ in $\gamma.\alpha$ with the
umbra $\gamma.\beta,$ we obtain the {\it composition} umbra of
$\alpha$ and $\gamma,$ i.e. $\gamma.\beta.\alpha.$ The
compositional inverse of an umbra $\alpha$ is the umbra
$\alpha^{<-1>}$ such that $\alpha^{<-1>}.\beta.\alpha \equiv \chi
\equiv \alpha.\beta.\alpha^{<-1>}.$ In particular it is
\begin{equation}
u^{<-1>}.\beta \equiv \chi \equiv u^{<-1>}.\beta,
\label{(chiinv)}
\end{equation}
where $u^{<-1>}$ denotes the compositional inverse of $u.$ Via the
umbral Lagrange inversion formula (cf. \cite{Dinardo}), the
moments of $u^{<-1>}$ are $E[(u^{<-1>})^n] =(-1)^{n-1}(n-1)!.$
Finally it is
\begin{equation}
\chi.\chi \equiv u^{<-1>} \qquad -\chi.(-\chi) \equiv (-u)^{<-1>}.
\label{(fattchi)}
\end{equation}
The disjoint sum of $\alpha$ and $\gamma$ is the umbra whose moments are the
sum of $n$-th moments $a_n$ and $g_n$ of $\alpha$ and $\gamma$
respectively, in symbols $(\alpha \dot{+} \gamma)^n \simeq \alpha^n + \gamma^n$
(cf. \cite{Dinardo1}). For instance it turns out
\begin{equation}
\chi.\alpha \dot{+} \chi.\gamma \equiv \chi.(\alpha+\gamma)
\label{(chiprop)}
\end{equation}
that is the well known additive property of cumulants. In the
following, we denote by  $\dot{+}_n \alpha$ the disjoint sum of
$n$ times the umbra $\alpha.$
\section{Umbral symmetric polynomials}
A partition of an integer $m$ is a sequence
$\lambda=(\lambda_1,\lambda_2,\ldots,\lambda_t),$ where $\lambda_i$
are weakly decreasing and $\sum_{i=1}^t \lambda_i=m.$ The integers
$\lambda_i$ are said {\it parts} of $\lambda.$ A different
notation is $\lambda=(1^{r_1},2^{r_2},\ldots),$ where $r_i$ is the
number of parts of $\lambda$ equal to $i.$ The {\it monomial
symmetric polynomials} in the variables
$\alpha_1,\alpha_2,\cdots,\alpha_n$ are $m_{\lambda}= \sum
\alpha_1^{\lambda_1} \cdots \alpha_t^{\lambda_t},$ where the sum
is over all distinct monomials having exponents
$\lambda_1,\ldots,\lambda_t.$ When $\lambda$ ranges in the set of
partition of the integer $m,$ $m_{\lambda}$ is a bases for the
algebra of the symmetric polynomials. There are different other
bases; we recall just those necessary in the following: the $r$-th
power sum symmetric polynomials $s_r=\sum_{i=1}^n \alpha^r_i;$ the
$k$-th elementary symmetric polynomials $e_k= \sum \alpha_{j_1}
\alpha_{j_2} \cdots \alpha_{j_k},$ where the sum is over $1 \leq
j_1<j_2<\cdots<j_k \leq n,$ and the $r$-th complete homogeneous
symmetric polynomials $h_i = \sum_{|\lambda|=i} m_{\lambda}.$
\par
When the umbrae $\alpha_1,\cdots,\alpha_n$ are uncorrelated and
similar to each other, these four classical bases of the symmetric
polynomial algebra can be represented by means of the umbral
polynomial $n.(\chi \alpha^i)$ and its moments. Propositions 1, 3
and 4 are stated under this hypothesis. For the power sum
symmetric polynomials it is
$$s_r \simeq n.\alpha^r \simeq n.(\chi \alpha^r).$$
Note that $s_r$ are umbrally equivalent to the moments of
$\dot{+}_n \alpha$ (cf. \cite{Dinardo1}).
\begin{proposition}[Umbral elementary polynomials]\label{el}
\begin{equation}
[n.(\chi\alpha)]^k \simeq \left\{ \begin{array}{cl}
k! e_k, &
k=1,2,\ldots,n, \\
0, & k=n+1,n+2,\ldots.
\end{array} \right.
\label{(ele)}
\end{equation}
\end{proposition}
\begin{proof}
For $k=1,\ldots,n$ the result follows applying the evaluation $E$
to the multinomial expansion of $[n.(\chi\alpha)]^k \simeq (\chi_1
\alpha_1 + \cdots + \chi_n \alpha_n)^k$ and observing that terms
having powers of $\chi$ greater than $1$ vanish. So just $k!$
monomials of the form $\chi_{j_1} \alpha_{j_1} \chi_{j_2}
\alpha_{j_2} \cdots \chi_{j_k} \alpha_{j_k}$ has an evaluation not
zero. Instead for $k=n+1,n+2,\ldots$ the result follows observing
that at least one power of $\chi$ greater than $1$ occurs in each
monomial of the multinomial expansion.
\end{proof}
\begin{proposition} It is
\begin{equation}
\chi.n.(\chi\alpha) \equiv u^{<-1>} \left( \dot{+}_n \alpha \right)
\quad n.(\chi\alpha) \equiv \beta. \left[u^{<-1>} \left( \dot{+}_n \alpha \right) \right]
\label{(powel)}
 \end{equation}
\end{proposition}
\begin{proof}
The first equivalence in (\ref{(powel)}) follows from (\ref{(chiprop)})
observing that
$$\chi.n.(\chi\alpha)\equiv \chi.(\chi_1\alpha_1+\cdot+\chi_n\alpha_n)
\equiv \dot{+}_n \chi.(\chi \alpha) \equiv \dot{+}_n u^{<-1>} \alpha.$$
The second one follows from the first, making the right dot-product
with $\beta$ and recalling (\ref{(chibeta)}).
\end{proof}
Equations (\ref{(powel)}) are the umbral version of the
well-known relations between power sum symmetric polynomials $s_r$
and elementary symmetric polynomials $e_k.$ Indeed $s_r$ is umbrally
equivalent to the moments of $(\dot{+}_n \alpha)$ and $e_k$
is umbrally equivalent to the moments of $n.(\chi \alpha).$ The
umbral expression of $m_{\lambda}$ requires the introduction of
augmented monomial symmetric polynomials $\tilde{m}_{\lambda}$.
Let $\lambda=(1^{r_1},2^{r_2},\ldots)$ be a partition of the
integer $m,$ such polynomials are defined by
\begin{equation}
\tilde{m}_{\lambda} = \sum_{1 \leq j_1 \ne \ldots
\ne j_{r_1} \ne j_{r_1+1} \ldots j_{r_1+r_2} \ldots \leq n} \alpha_{j_1} \cdots \alpha_{j_{r_1}}
\alpha^2_{j_{r_1+1}} \cdots \alpha^2_{j_{r_1+r_2}}  \cdots.
\label{(aug)}
\end{equation}
The proof of the next proposition follows the same patterns of
Proposition \ref{el}.
\begin{proposition} It is
$$\tilde{m}_{\lambda} \simeq
[n.(\chi\alpha)]^{r_1}[n.(\chi \alpha^2)]^{r_2} \cdots.$$
\end{proposition}
The next corollary follows recalling that $m_{\lambda} =\tilde{m}_{\lambda}/[r_1! r_2!
\cdots].$
\begin{corollary} [Umbral monomial polynomials] \label{cr}
$$m_{\lambda} \simeq
\frac{[n.(\chi\alpha)]^{r_1}}{r_1!} \frac{[n.(\chi\alpha^2)]^{
r_2}}{r_2!} \cdots.$$
\end{corollary}
\begin{proposition}[Umbral complete polynomials]
\begin{equation}
[-n.(-\chi\alpha)]^m \simeq m! h_m \qquad m=1,2,\ldots.
\label{(compl)}
\end{equation}
\end{proposition}
\begin{proof}
It is
\begin{eqnarray*}
[-n.(-\chi\alpha)]^m & \simeq & [-1.(-\chi_1 \alpha_1) + \cdots +
-1.(-\chi_n \alpha_n)]^m \\
& \simeq & m! \sum_{|\lambda|=m} [-1.(-\chi{'})]^{.r_1} ([-1.(-\chi^{''})]^2)^{.r_2}
\cdots \frac{\tilde{m}_{\lambda}}{(1!)^{r_1} r_1!(2!)^{r_2}r_2! \cdots}
\end{eqnarray*}
and the result follows from Corollary \ref{cr}, being $([-1.(-\chi)]^i)^{.r_i}\simeq (i!)^{r_i}.$
\end{proof}
Note that equivalences (\ref{(ele)}) and (\ref{(compl)}) are the
umbral version of the well-known identities:
$$\sum_{k} e_k t^k = \prod_{i=1}^n (1+\alpha_i t) \qquad \sum_{k} h_k t^k =
\frac{1}{\prod_{i=1}^n (1-\alpha_i t)}.$$
\begin{proposition} It is
\begin{equation}
-\chi.n.(-\chi\alpha) \equiv (-u)^{<-1>} \left( \dot{+}_n \alpha \right)
\quad -n.(-\chi\alpha) \equiv \beta. \left[(-u)^{<-1>} \left( \dot{+}_n \alpha \right) \right]
\label{(powel2)}
 \end{equation}
\end{proposition}
\begin{proof}
The results follow from  (\ref{(powel)}) replacing the umbra $\chi$
with $-\chi$ and recalling that $u^{<-1>} \equiv \chi.\chi$ must
be replaced with $(-u)^{<-1>} \equiv -\chi.(-\chi).$
\end{proof}
Equations (\ref{(powel2)}) are the umbral version of the
well-known relations between power sum symmetric polynomials $s_r$
and complete symmetric polynomials $h_k.$ Indeed $s_r$ is umbrally
equivalent to the moments of $(\dot{+}_n \alpha)$ and $h_k$ is
umbrally equivalent to the moments of $-n.(-\chi \alpha).$
\section{$U$-statistics}
Let $X_1, X_2, \ldots, X_n$ be $n$ independent r.v.'s.
A statistic of the form
$$U=\frac{1}{(n)_k}\sum \Phi(X_{j_1}, X_{j_2}, \ldots, X_{j_k})$$
where the sum ranges in the set of all permutations
$(j_1,j_2,\ldots,j_k)$ of $k$ integers, $1 \leq j_i \leq n,$ is
called $U-$statistic \cite{Hoeffding}. If $X_1,X_2, \ldots, X_n$
have the same cumulative distribution function $F(x),$ $U$ is an
unbiased estimator of the population character $\theta(F)=\int
\cdots \int \Phi(x_1,\cdots,x_k) dF(x_1)\cdots dF(x_k).$ In this
case, the function $\Phi$ may be assumed to be a symmetric
function of its arguments. Often, in the applications, $\Phi$ is a
polynomial in $X_i$'s so that the $U-$statistic is a symmetric
polynomial. By virtue of the fundamental theorem on symmetric
polynomials, the $U-$statistic can be considered a polynomial in
the elementary symmetric polynomials. The following theorem is an
umbral reformulation of the above statement.
\begin{theorem}[$U-$statistic]
If $\lambda=(1^{r_1},2^{r_2},\ldots,)$ is a partition of
the integer $m \leq n$ then
\begin{equation}
(\alpha_{j_1})^{.r_1} (\alpha_{j_2}^2)^{.r_2} \cdots  \simeq
\frac{1}{(n)_k} [n.(\chi\alpha)]^{r_1}[n.(\chi\alpha^2)]^{r_2}\cdots
\label{(ustat)}
\end{equation}
where $j_i \in \{1,2,\ldots,n\}$ and $\sum_{j} r_j=k.$
\end{theorem}
\begin{proof}
The result follows observing that $[n.(\chi\alpha^i)]^{r_i} \simeq
(n.\chi)^{r_i} (\alpha^{i})^{.r_i}$ and $(n.\chi)^k \simeq (n)_k.$
\end{proof}
The formula (\ref{(ustat)}) states how to estimate moment products
by using only $n$ information  drawn out the population. Then the
symmetric polynomial on the right side of (\ref{(ustat)}) is the
$U-${\it statistic} of the uncorrelated and similar umbrae
$\alpha_1,\alpha_2,\ldots,\alpha_n$ associated to
$(\alpha_{j_1})^{.r_1} (\alpha_{j_2}^2)^{.r_2} \cdots.$
\begin{example}{\it Moment powers.}
{\rm Set $r_1=2$ and $k=2,$ from (\ref{(ustat)}) the $U-$statistic
associated to $\alpha^{.2} \simeq a_1^2$ is
$$\alpha^{.2} \simeq \frac{1}{n(n-1)} [n.(\chi \alpha)]^2 \simeq \frac{1}{n(n-1)} \sum_{i\ne j}
\alpha_i \alpha_j, \qquad n \geq 2.$$}
\end{example}
\begin{example}{\it $h$-statistics.}
{\rm As it has been shown in \cite{Dinardo1}, it is
$$(\alpha^{a_1})^r \simeq \sum_{k=0}^r \left( \begin{array}{c}
r \\
k
\end{array} \right) (-1)^k \alpha^{.k} (\alpha^{'})^{r-k}
\simeq \sum_{k=0}^{r-2} \frac{(-1)^k}{r-k} \sum_{k=0}^r \alpha^{.k} (\alpha^{'})^{r-k}$$
where $\alpha^{a_1}$ is the central umbra of $\alpha$ about $a_1=E[\alpha]$
and $\alpha^{'} \equiv \alpha$ is an umbra uncorrelated with $\alpha.$
Replacing the product $\alpha^{.k} (\alpha^{'})^{r-k}$
with the corresponding $U-$statistic (\ref{(ustat)}), it results
\begin{equation}
(\alpha^{a_1})^r  \simeq  \sum_{k=0}^r \left( \begin{array}{c}
r \\
k
\end{array} \right) \frac{(-1)^k}{(n)_{k+1}}  [n.(\chi\alpha)]^{k} n.(\chi
\alpha^{r-k}).
\label{(hstat1)}
\end{equation}
When the products $[n.(\chi\alpha)]^{k} n.(\chi \alpha^{r-k})$ are
expressed in terms of power sum symmetric polynomials, we have the
so called $h$-statistics. On this matter, we will give more
details in the next section.}
\end{example}
\section{$k$-statistics}
The $n$-th $k$-statistic $k_n$ is the unique symmetric unbiased
estimator of the $n$-th cumulant $\kappa_n$ of a given statistical
distribution, i.e. $E[k_n]=\kappa_n.$ The $k$-statistics can be
expressed in terms of the sums of the $r$-th powers of the data
points. In this section we give an umbral syntax that provides a
general computational method to generate such expressions.
To this aim, we digress to introduce the exponential Bell umbral polynomials.
\par
The most widespread expression of incomplete exponential Bell
polynomials is referred to partition of an integer. Of course, it
is also possible to express such polynomials referring them to
partition of a set. Here we follow this last point of view. We
denote by $\Pi_{i,k}$ the set of all partitions of the set
$[i]=\{1,2,\ldots,i\}$ in $k$ blocks. Let $\pi=\{A_1,A_2,
\ldots,A_k\}$ be an element of $\Pi_{i,k}$. Then it is
$B_{i,k}(a_1,a_2,\ldots) = \sum_{\pi \in \Pi_{i,k}} a_{n_1} \,
a_{n_2} \, \cdots \, a_{n_k}$ where $|A_j|=n_j, j=1,2,\ldots,k,$ as we will suppose
from now on.  Let us consider the following symmetric umbral
polynomial:
\begin{equation}
{\cal B}_{i,k}(\alpha_1,\alpha_2,\cdots, \alpha_k) = \sum_{\pi \in
\Pi_{i,k}} \alpha_1^{n_1} \, \alpha_2^{n_2} \, \cdots \,
\alpha_k^{n_k}, \label{(bellexp2u)}
\end{equation}
where $\alpha_1,\alpha_2,\cdots, \alpha_k$ are uncorrelated umbrae
similar to $\alpha$ and $a_{n_j}=E[\alpha^{n_j}], j=1,2,\ldots,k
.$ Obviously it is
$E[{\cal B}_{i,k}]=B_{i,k},$ so that any expression containing the
polynomials $B_{i,k}$ could be replaced with an umbrally
equivalent expression containing the polynomials ${\cal B}_{i,k}.$
The polynomials ${\cal B}_{i,k}$ will be called (incomplete)
umbral exponential Bell polynomials.  The combinatorics underlain
the polynomial ${\cal B}_{i,k}$ is the following: the set $[i]$ is
partitioned in $k$ blocks, to each of them one associates the
umbra $\alpha^{n_j}$ obtained firstly replacing the elements in
the $j$-th block with the umbra $\alpha$ and then labelling all blocks so
that powers belonging to different blocks result uncorrelated.
Replacing in (\ref{(bellexp2u)}) the products $\alpha_1^{n_1} \,
\alpha_2^{n_2} \, \cdots \, \alpha_k^{n_k}$ with the umbrally
corresponding $U-$statistic (\ref{(ustat)}), we get for $i \leq n$
\begin{equation}
{\cal B}_{i,k}(\alpha_1,\alpha_2,\cdots, \alpha_k) \simeq
\sum_{\pi \in \Pi_{i,k}}  \frac{1}{(n)_k} n.(\chi \, \alpha^{n_1})
\, n.(\chi \, \alpha^{n_2})  \, \cdots \, n.(\chi \,
\alpha^{n_k}), \label{(bellustat)}
\end{equation}
by which we are able to give the umbral $k$-statistics. Indeed, the
$\alpha-$cumulant umbra $\kappa_{\alpha}$ is similar to
$\chi.\alpha,$ so that, being $(\chi)_k \simeq (u^{<-1>})^k,$
from (\ref{(ombdot)}) and (\ref{(bellustat)}) it is
\begin{equation}
(\chi.\alpha)^i  \simeq  \sum_{k=1}^{i} (-1)^{k-1}
\frac{(k-1)!}{(n)_k} \sum_{\pi \in \Pi_{i,k}} n.(\chi \,
\alpha^{n_1}) \, n.(\chi \, \alpha^{n_2})  \, \cdots \, n.(\chi \,
\alpha^{n_k}). \label{(kstat)}
\end{equation}
Since $n.(\chi \alpha^{n_i})$ is umbrally equivalent to a
symmetric power sum polynomial, equation (\ref{(kstat)}) gives
the moments of the $\alpha-$cumulant umbra in terms of power sum polynomials, i.e. the
umbral form of $k$-statistics. Note that the symmetric power sum
polynomials in (\ref{(kstat)}) are correlated. So, in order to
make formula (\ref{(kstat)}) effective, we need a device by
which to evaluate the product $ n.(\chi \,
\alpha^{n_1}) \, n.(\chi \, \alpha^{n_2})  \, \cdots \, n.(\chi \,
\alpha^{n_k})$. To this aim, by using the umbral exponential Bell polynomials
(\ref{(bellexp2u)}), the moments of the umbra $n.(\chi \alpha)$ can be evaluate from
(\ref{(ombdot)}) recalling $n.\alpha^{n_i} \simeq (\dot{+}_n
\alpha)^{n_i}$ and the second equivalence in (\ref{(powel)})
\begin{eqnarray}
[n.(\chi \alpha)]^i \simeq \sum_{k=1}^{i} \sum_{\pi \in \Pi_{i,k}}
(u^{'<-1>})^{n_1}(u^{''<-1>})^{n_2} \cdots (u^{'''<-1>})^{n_k}
[n. {\alpha^{'}}^{n_1}] \, [n.{\alpha^{''}}^{n_2}] \cdots
[n.{\alpha^{'''}}^{n_k}].
\label{(momdotbeta1)}
\end{eqnarray}
The previous umbral equivalence is suitable to be generalized
to the product of umbral polynomials $[n.(\chi p_1)] [n.(\chi p_2)] \cdots [n.(\chi p_i)]$
with no disjoint support. Indeed, it results
\begin{equation}
[n.(\chi p_1)] \cdots [n.(\chi p_i)] \simeq \sum_{k=1}^{i} \sum_{\pi \in \Pi_{i,k}}
(u^{'<-1>})^{n_1} \cdots (u^{''<-1>})^{n_k}
[n. P^{'}_{A_1}] \, \cdots \, [n.P^{''}_{A_k}]
\label{(dotprod)}
\end{equation}
where $P_{A_j}=\prod_{t=1}^{n_j}p_{j_t}$ and $p_{j_t}$ are the
polynomials indexed by the elements of the block $A_j,$ as we
will suppose from now on.  Equivalence (\ref{(dotprod)}) is the
device required to generate $k$-statistics.
Indeed, setting $i=k$ and $p_t=\alpha^{n_t}$ for
$t=1,2,\cdots,i,$ it results:
\begin{equation}
[n.(\chi \alpha^{n_1})] \cdots [n.(\chi \alpha^{n_k})]
\simeq \sum_{j=1}^{k} \sum_{\pi \in \Pi_{k,j}}
(u^{'<-1>})^{n_1} \cdots (u^{''<-1>})^{n_j}
[n. {\alpha^{'}}^{m_1}]  \cdots [n.{\alpha^{''}}^{m_j}]
\label{(cum5)}
\end{equation}
where $m_j = \sum_{t=1}^{n_j}n_{j_t}$ and $n_{j_t}$ are indexed
by the elements of the block $A_j.$  Note that the power sum
polynomials on the right side of (\ref{(cum5)}) are now uncorrelated so such equivalence
gives augmented monomial symmetric polynomials in
terms of power sum polynomials, translating Kendall and Stuart tables
read downwards \cite{Stuart-Ord}. Instead the following formula leads to a
symbolic translation of the Kendall and Stuart tables read across
(i.e. power sum polynomials in terms of augmented monomial ones):
\begin{equation}
(n.p_1) \cdots(n.p_i) \simeq  \sum_{k=1}^i (n)_k
\sum_{\pi \in \Pi_{i,k}} P^{'}_{A_1} \cdots P^{''}_{A_k} \simeq
\sum_{k=1}^{i} \sum_{\pi \in \Pi_{i,k}} n.(\chi P_{A_1}) \cdots
n.(\chi P_{A_k}).
\label{(momdotchi1)}
\end{equation}
The first equivalence in (\ref{(momdotchi1)}) is
obtained from (\ref{(momdot)}) trough analogous considerations
used to state (\ref{(dotprod)}); the second equivalence comes from
(\ref{(ustat)}) replacing $\alpha_i^{j_i}$ with the umbral polynomial $p_i,$
i.e.
$$ P^{'}_{A_1} P^{''}_{A_2} \cdots
P^{'''}_{A_k} \simeq \frac{1}{(n)_k}n.(\chi P_{A_1}) n.(\chi
P_{A_2}) \cdots n.(\chi P_{A_k}).$$
Set in (\ref{(momdotchi1)}) $p_t=\alpha^{n_t}$ for $t=1,2,...,i,$ we have
\begin{equation}
(n.\alpha^{n_1}) \cdots(n.\alpha^{n_i}) \simeq  \sum_{k=1}^{i}
\sum_{\pi \in \Pi_{i,k}} (n.\chi \alpha^{m_1})(n.\chi \alpha^{m_2}) \cdots (n.\chi
{\alpha}^{m_k}) \label{(momdotchi2)}
\end{equation}
where $m_j = \sum_{t=1}^{n_j}n_{j_t}$ and $n_{j_t}$ are indexed
by the elements of the block $A_j.$
\begin{example}{\it $h$-statistics.}
{\rm In (\ref{(cum5)}) set $n_1=n-k$ and $n_2=\cdots=n_{k+1}=1;$
from (\ref{(hstat1)}) we get the umbral
expression of $h$-statistics.}
\end{example}
\begin{example}{\it Joint cumulants.}
{\rm Let  $p_1,p_2,\ldots,p_i$ be umbral polynomials. Replacing $n$ with $\chi$
in the first equivalence in (\ref{(momdotchi1)}), we have
\begin{equation}
(\chi.p_1)(\chi.p_2)\cdots(\chi.p_i)  \simeq   \sum_{k=1}^i (\chi)_k
\sum_{\pi \in \Pi_{i,k}} P^{'}_{A_1} \cdots P^{''}_{A_k}. \label{(cum3)}
\end{equation}
When the umbral polynomials $p_i$ are interpreted as r.v.'s, equivalence (\ref{(cum3)})
gives their joint cumulants. So we will call $(\chi.p_1)(\chi.p_2)\cdots(\chi.p_i)$
the joint cumulant of $p_1,\ldots,p_i.$ Note that, setting $p_t=\alpha$
for $t=1,2,\ldots,i,$ one has the $i$-th ordinary cumulant
$(\chi.\alpha)^i.$ Through this equivalence it results
$\chi.(p_1 + \cdots + p_i) \equiv \chi.p_1 + \cdots + \chi.p_i.$
Now suppose to split the set $\{p_1,p_2,\ldots,p_i\}$ in two
subsets $\{p_{j_1},\ldots,p_{j_t}\}$ and $\{p_{k_1},\ldots,p_{k_s}\}$
with $s+t=i,$ such that polynomials belonging to different subsets
are uncorrelated. Then it is
\begin{equation}
(\chi.p_1)(\chi.p_2)\cdots(\chi.p_i) \simeq 0.
\label{(cumconj)}
\end{equation}
Indeed, setting $P=\sum_{l=1}^t p_{j_l}$ and $Q=\sum_{l=1}^{s}
p_{k_l},$ such polynomials are uncorrelated, so that $\chi.(P + Q)
\equiv \chi.P \dot{+} \chi.Q$ from (\ref{(chiprop)}). Equivalence
(\ref{(cumconj)}) follows observing that, due to the disjoint sum,
products involving powers of $\chi.P$ and $\chi.Q$ vanish. When
the umbral polynomials $p_i$ are interpreted as r.v.'s,
equivalence (\ref{(cumconj)}) states the following well-known
result:  if some of the r.v.'s are uncorrelated of all others,
then their joint cumulant is zero.}
\end{example}
\begin{example}{\it Multivariate $k$-statistics.}
{\rm Equivalence (\ref{(momdotchi1)}) allows a compact expression
of multivariate $k$-statistics. Replacing $n$ with $\chi$ in its
second equivalence, we construct the $U-$statistic of the joint
cumulant
\begin{equation}
(\chi.p_1)(\chi.p_2)\cdots(\chi.p_i)  \simeq
\sum_{k=1}^{i} \frac{(\chi)_k}{(n)_k} \sum_{\pi \in \Pi_{i,k}}
n.(\chi P_{A_1}) n.(\chi P_{A_2}) \cdots n.(\chi P_{A_k}) \label{(cum4)}.
\end{equation}
Again, in the product on the right side of (\ref{(cum4)}) the umbral polynomials
are correlated. In order to make effective the computation, it is necessary to rewrite (\ref{(cum4)}) by using
equivalence (\ref{(dotprod)}) with $P_{A_t}$ instead of $p_t.$
For instance, in order to express $k_{21},$ set in (\ref{(cum4)}) $i=3$ and
$p_1=p_2=\alpha_1, p_3=\alpha_2.$ It results
\begin{eqnarray}
(\chi.\alpha_1)^2(\chi.\alpha_2) & \simeq & \frac{\chi}{n} \,
n.(\chi \alpha^2_1 \alpha_2) + \frac{(\chi)_2}{(n)_2} \, \left\{ 2
\, n.(\chi \alpha_1) \, n.(\chi \alpha_1 \alpha_2)
+ \, n.(\chi \alpha^2_1)\, n.(\chi \alpha_2) \right\} \nonumber \\
& + & \frac{(\chi)_3}{(n)_3} \left\{ [n.(\chi \alpha_1)]^2 \, n.(\chi
\alpha_2)\right\} \label{(k21)}
\end{eqnarray}
Set $s_{p,q} \simeq n.(\alpha_1^p \, \alpha_2^q).$
It is
\begin{eqnarray}
n.(\chi \alpha_1) \, n.(\chi \alpha_1 \alpha_2) & \simeq & (u^{<-1>})^2 n.(\alpha_1^2
\alpha_2) + (u^{<-1>})^{.2} n.\alpha'_1 \, n.(\alpha_1 \alpha_2)
\simeq - s_{2,1} + s_{1,0}\, s_{1,1} \label{(1)}  \\
n.(\chi \alpha^2_1) \, n.(\chi \alpha_2) & \simeq & (u^{<-1>})^2 n.(\alpha_1^2
\alpha_2) + (u^{<-1>})^{.2} n.{\alpha'}^{2}_1 \, n.\alpha_2
\simeq - s_{2,1} + s_{2,0} \, s_{0,1} \label{(2)} \\
\{n.(\chi \alpha_1)\}^2 n.(\chi \alpha_2) & \simeq & (u^{<-1>})^3 n.(\alpha_1^2
\alpha_2) + (u^{<-1>})^{.3} n.{\alpha'}_1 \, n.\alpha_1 \, n.\alpha_2 \nonumber \\
& + & {u^{'}}^{<-1>} \, (u^{<-1>})^{2} [n.\alpha^{2}_1 \, n.\alpha_2 + 2
\, n.\alpha^{'}_1 \, n.(\alpha_1 \alpha_2)] \nonumber \\
& \simeq &  2 \, s_{2,1} - s_{2,0} \, s_{0,1} -2 \, s_{1,0} \, s_{1,1} + s^2_{1,0}\,  s_{0,1}
\label{(3)}
\end{eqnarray}
Equivalence (\ref{(1)}) comes from (\ref{(dotprod)}) setting $i=2,$
$p_1=\alpha_1$ and $p_2=\alpha_1 \alpha_2,$  equivalence (\ref{(2)})
comes from (\ref{(dotprod)}) setting $i=2,$ $p_1=\alpha^2_1$ and $p_2=\alpha_2,$
equivalence (\ref{(3)}) comes from (\ref{(dotprod)}) setting $i=3,$ $p_1=p_2=\alpha_1$
and $p_3=\alpha_2.$ Substituting the above equivalences in
(\ref{(k21)}) and rearranging the terms, we have the expression of $k_{21},$
$$k_{21} \simeq (\chi.\alpha_1)^2(\chi.\alpha_2) \simeq \frac{1}{(n)_3} \left[n^2 s_{2,1} - 2 n \, s_{1,0}\, s_{1,1} - n \, s_{2,0}
\, s_{0,1} + 2\, s^2_{1,0}\, s_{0,1} \right].$$ The expression of
generalized $k-$ statistics (as well as the multivariate ones) in
terms of power sums comes from (\ref{(cum3)}) replacing some of
the umbrae $\chi$ with uncorrelated ones and then constructing the
corresponding $U-$statistics.}
\end{example}


\begin{thebibliography}{99}
\setlength{\baselineskip}{10pt}
%
\bibitem{Dinardo}
E. Di Nardo, D. Senato, `Umbral nature of the Poisson random variables',
In Algebraic combinatorics and computer science (eds. Crapo H. and Senato
D.), Springer Italia (2001), 245--266.
%
\bibitem{Dinardo1}
E. Di Nardo, D. Senato, `An umbral setting for cumulants and factorial moments',
{\it Europ. Jour. Combinatorics }, (2005) to appear.
%
\bibitem{Doubilet}
P. Doubilet, `On the foundations of combinatorial theory VII:
Symmetric functions through the theory of distribution and
occupancy',  {\it Stud. Appl. Math.}  {\bf 11}, (1972) 377--396.
%
\bibitem{Fisher}
R.A. Fisher, `Moments and product moments of sampling
distributions', {\it Proc. London Math. Soc.} (2) {\bf 30}, (1929)
199--238.
%
\bibitem{Hoeffding}
W. Hoeffding, `A class of statistics with asymptotically normal distribution',
{\it Ann. Math. Stat.} {\bf 19}, (1948) 293--325.
%
\bibitem{SIAM}
G.-C. Rota, B.D. Taylor, `The classical umbral calculus',
{\it SIAM J. Math. Anal.} {\bf 25}, (1994) 694--711.
%
\bibitem{Speed2}
T. P. Speed, `Cumulants and partition lattices II: Generalized
$k$-statistics', {\it J. Aust. Math. Soc.}, Ser. A {\bf 40}, (1986) 34--53.
%
\bibitem{Stuart-Ord}
A. Stuart, J.K. Ord, {\it Kendall's Advanced Theory of
Statistics}, Vol. 1, Charles Griffin and Company Limited, London,
(1987).
%
\bibitem{Taylor1}
B.D. Taylor, `Difference equations via the classical umbral
calculus' In Mathematical Essays in Honor of Gian-Carlo Rota (eds. Sagan
et al.), Birkhauser Boston (1998), 397--411.
\end{thebibliography}
\end{document}